\documentstyle[12pt]{article}
\title{  CLASSICAL YANG-BAXTER EQUATION AND
 LOW DIMENSIONAL TRIANGULAR LIE BIALGEBRAS
 \thanks {This work is supported by National Science Foundation}}
\author{Shouchuan Zhang \\ Department  of Mathematics, Nanjing  University,
210008 \\ P.R.China}
\date{}
\begin{document}
\newtheorem{Theorem}{\quad Theorem}[section]
\newtheorem{Proposition}[Theorem]{\quad Proposition}
\newtheorem{Definition}[Theorem]{\quad Definition}
\newtheorem{Corollary}[Theorem]{\quad Corollary}
\newtheorem{Lemma}[Theorem]{\quad Lemma}
\newtheorem{Example}[Theorem]{\quad Example}
\maketitle
\begin {abstract}
 All  solutions of constant classical Yang-Baxter
 equation (CYBE)  in Lie algebra $L$ with dim $L \le 3$  are obtained
  and
 the sufficient and  necessary conditions which  $(L, \hbox {[ \ ]},
 \Delta _r, r)$ is a coboundary
 (or triangular ) Lie bialgebra  are given.
The strongly symmetric elements in $L\otimes L $ are found and they all are
 solutions of CYBE in $L$ with $dim \ \ L \le 3$.
\end {abstract}
\noindent
\addtocounter{section}{-1}

\section{ Introduction}

The concept and structures of Lie coalgebras were  introduced
and  studied by   W.Michaelis
in \cite {M80}
\cite {M85}.
 V.G.Drinfel'd and   A.A.Belavin in    \cite {D1}
\cite {BD}  introduced the notion  of
triangular, coboundary  Lie bialgebra $L$ associated to a solution
$r \in L \otimes L$ of the CYBE and gave a classification of  solutions of CYBE with
parameter for  simple Lie algebras. W.Michaelis in \cite {M} obtained the
structure of
a triangular, coboundary Lie bialgebra on any Lie algebra containing
linearly independent elements $a$ and $b$ satisfying $[a,b]= \alpha b$
for some non-zero $\alpha \in k$ by setting
$r= a\otimes b - b \otimes a $.

The Yang-Baxter equation first came up in a paper by Yang as
factorition condition of the scattering S-matrix in
the many-body  problem in one dimension and in work of Baxter on
exactly solvable models in statistical mechanics.
It has been playing an important role in mathematics and physics
( see \cite {BD} , \cite {YG} ).
Attempts to find solutions of The Yang-Baxter equation in a systematic
way have let to the theory of quantum groups.
The Yang-Baxter equation is of many forms.
The classical Yang-Baxter equation is one.

 In many applications one need to know the solutions of classical Yang-Baxter
 equation and know if a Lie algebra is a coboundary Lie bialgebra or
 a triangular Lie bialgebra.  A systematic study of low dimensional Lie
 algebras,
  specially, of those Lie algebras that play a role in physics
 (as e.g. $sl (2, {\bf C})$, or the Heisenberg algebra), is very useful.

 In this paper, we obtain all  solutions of constant classical Yang-Baxter
 equation (CYBE)  in Lie algebra $L$ with dim $L \le 3$
  and  give
 the sufficient and  necessary conditions which  $(L, \hbox {[ \ ]},
 \Delta _r, r)$ is a coboundary
 (or triangular ) Lie bialgebra.
We find the strongly symmetric elements in $L\otimes L $ and show
they all are
 solutions of CYBE in $L$ with $dim \ \ L \le 3$.
  Using these conclusions, we study the
 Lie algebra  $sl(2)$.

In order to make the paper somewhat self-contained, we begin by recalling
the definition of Lie bialgebra (all of the following definitions in
the introduction are the same as \cite [P368-371] {M}).

Let $k$ be a field with char $k \not= 2 $ and $L$ be a Lie algebra over $k$.
If $r = \sum a_i \otimes b_i  \in L \otimes L $ and $x\in L$, we define
 \begin{eqnarray*}                              
\left[ r^{12}, r^{13} \right] &:=& \sum_{i, j} [a_i, a_j] \otimes b_i
\otimes b_j \\
\left[ r^{12}, r^{23} \right] &:=& \sum_{i, j} a_i\otimes [b_i,a_j]
 \otimes  b_j \\
\left[ r^{13}, r^{23} \right] &:=& \sum_{i, j} a_i\otimes a_j \otimes [b_i, b_j] \\
       x \cdot r &:=& \sum _{i} [x,a_i] \otimes b_i + a_i \otimes [x,b_i] \\
      \Delta _r (x) &:=& x \cdot r
\end{eqnarray*}

and call $$     [r^{12}, r^{13}]+  [r^{12}, r^{23}]+ [r^{13}, r^{23}]=0$$
the classical Yang-Baxter equation (CYBE).

Let $$ \tau: L \otimes L \longrightarrow L \otimes L  $$
denote the natural twist map ( defined by $x \otimes y \mapsto y \otimes x$ ),
 and let $$
 \xi : L \otimes  L \otimes L \longrightarrow L \otimes  L \otimes L  $$
 the map defined by
 $ x \otimes y \otimes z \mapsto  y \otimes z \otimes x$
 for any $x, y, z \in L$.

a vector space $L$ is called a Lie coalgebra, if there exists a
linear map $$\Delta : L \longrightarrow L \otimes L $$
such that

(i)  $Im \Delta \subseteq Im (1- \tau )$ and

(ii) $(1 + \xi + \xi ^2)(1 \otimes \Delta ) \Delta =0$

a vector space $(L, [\hbox { \ }], \Delta )$  is called a Lie bialgebra
if

 (i)  $(L,[\hbox { \ }])$   is a Lie algebra,

(ii)  $(L, \Delta ) $  is a Lie coalgebra,

(iii) for all $x, y \in L$
$$ \Delta [x,y]= x \cdot \Delta (y) - y \cdot \Delta (x), $$
 where, for all $x, a_i, b_i  \in L$,

$(L,[\hbox { \ }], \Delta, r )$   is called a coboundary Lie
 bialgebra, if  $(L,[\hbox { \ }], \Delta )$  is  a  Lie
bialgebra and  $r \in Im (1-\tau ) \subseteq L \otimes L $
such that
$$\Delta (x) = x \cdot r$$  for all $x\in L$.
a coboundary Lie bialgebra $(L,  [\hbox { \ }], \Delta, r )$ is called
triangular, if $r$ is a solution of CYBE.

\section {The solutions of CYBE }
In this section, we find the general solution of CYBE for Lie algebra $L$
with $dim \ \ L \le 3.$

Jacobson gave  a classification of Lie algebras with their dimension
$dim L \le 3$ in \cite [ P 11-14] {J}. We now write all of their operation
as follows:

(I)  If $L$ is an abelian Lie algebra,
then its operation is trivial;

(II) If  dim $ L $= dim $L'=3$, then there exist a basis
 $\{ e_1, e_2, e_3 \}$ and $\alpha, \beta \in k$ with $\alpha \beta \not=0$
such that
$[e_1, e_2]= e_3, [e_2, e_3]= \alpha  e_1,  [e_3, e_1]=\beta e_2$;

(III) If  dim $ L =3$ and $L' \subseteq $ the center of $L$ with
dim $L'=1$, then there exist a basis
 $\{ e_1, e_2, e_3 \}$ and $\alpha, \beta \in k$ with $\alpha= \beta =0$
such that
$[e_1, e_2]= e_3, [e_2, e_3]= \alpha  e_1,  [e_3, e_1]=\beta e_2;$

(IV) If  $k$ is an algebraically closed field and dim $ L =3$
 with dim $L'=2$, then there exist a basis
 $\{ e_1, e_2, e_3 \}$ and $\beta, \delta \in k$ with $\delta \not=0$
such that
$[e_1, e_2]= 0, [e_1, e_3]= e_1 + \beta  e_2,  [e_2, e_3]=\delta e_2;$

(V) If  dim $ L =3$ and $L' \not\subseteq$ the  center of $L$ with
dim $L' =1$,
 then there exist a basis
 $\{ e_1, e_2, e_3 \}$ and $\beta, \delta \in k$ with $\beta= \delta=0$
such that
$[e_1, e_2]= 0, [e_1, e_3]=  e_1 + \beta e_2,  [e_2, e_3]=\delta e_2;$

(VI) If  dim $ L=2$ with dim $L'=1$, then there exists a basis
 $\{ e,f \}$
such that
$[e,f]= e.$

We like to  point out the actual meaning of the different algebras, namely
(II) includes  $sl (2, {\bf C})$ and its real form $su (2)$, (III)
is just the three-dimensional Heisenberg algebra, (IV) includes the $(1+1) $
Poincare algebra, (V) is the (one-dimensional) central extension of the
non-abelian two-dimensional Lie algebra and (VI) is the non-abelian two-
dimensional Lie algebra itself.

For case (I), the two treated problems are trivial.

\begin {Definition} \label {1.1}
  Let  $\{ e_1, e_2, \cdots, e_n \}$ be a basis of vector space $V$
and     $r = \sum_{i,j= 1 }^{n} k_{ij}(e_{i} \otimes e_{j}) \in V \otimes V,$
where   $k_{ij} \in k,$ for  $i, j =1, 2, \cdots, n.$

(i) If  $$k_{ij}=k_{ji}  \hbox { \ \ \ } k_{ij} k_{lm} = k_{il}k_{jm}$$
for $i, j, l, m = 1, 2, \cdots, n,$   then  $r$ is called
strongly symmetric to the basis
 $\{ e_1, e_2, \cdots, e_n \}$;

(ii) If   $k_{ij} =- k_{ji}$, for $i,j = 1, 2, \cdots, n$,
 then  $r$ is called
 skew symmetric.

(iii) Let $dim { \ } V =3$ and
 $k_{11}=x, k_{22}=y,  k_{33}=z, k_{12}=p,  k_{21} =q, k_{13}=s,
 k_{31} =t, k_{23}=u,  k_{32}=v, \alpha, \beta \in k.$
If $p=-q,s=-t, u=-v,  x = \alpha z, y= \beta z$ and
$$\alpha \beta z^2 +   \beta s^2 + \alpha u^2 + p^2 =0,$$
then  $r$ is called
 $\alpha , \beta $-skew symmetric to the
basis  $\{ e_1, e_2, e_3 \}.$
\end {Definition}
Obviously, skew symmetry does not depend on the particular choice of basis
of $V$ . We need to know if
strong  symmetry  depend on the particular choice of basis of $V$.

\begin {Lemma}                                            
\label {1.2}
(I)  If $V$ is a finite-dimensional vector space,
then the strong symmetry  does not depend on the particular choice of
basis of $V;$

(II)  Let  $\{ e_1, e_2, e_3 \}$ be a basis of vector space $V,$
  $ p,q, s,t, u, v, x, y, z \in k $  and
    $r=  p(e_1 \otimes e_2)+q (e_2 \otimes e_1)
+ s (e_1 \otimes e_3)+t (e_3 \otimes e_1)
+u (e_2 \otimes e_3) +v (e_3 \otimes e_2) +x (e_1 \otimes e_1)
+y (e_2 \otimes e_2) + z (e_3 \otimes e_3).$
           
Then   $r$ is strongly symmetric iff
 $$ xy =p^2, xz = s^2,  yz=u^2, xu=sp, ys=pu, zp=su,
p=q, s=t, u=v$$ iff  $$xy = p^2, xz = s^2,  yz=u^2, xu=sp,
p=q, s=t, u=v$$ iff
 $$ xy =p^2, xz = s^2,  yz=u^2, ys=up,
p=q, s=t, u=v$$ iff
 $$ xy =p^2, xz = s^2,  yz=u^2, zp=su,
p=q, s=t, u=v$$ iff

 $r=   z^{-1}su  (e_1 \otimes e_2)+ z^{-1}su   (e_2 \otimes e_1)
+ s (e_1 \otimes e_3)+ s (e_3 \otimes e_1)
+u (e_2 \otimes e_3) +u (e_3 \otimes e_2)
+z^{-1}s^2  (e_1 \otimes e_1)
+ z^{-1} u^2 (e_2 \otimes e_2) + z (e_3 \otimes e_3)$
                 with $z \not= 0;$   or,

  $r= p(e_1 \otimes e_2)+ p (e_2 \otimes e_1)
+x (e_1 \otimes e_1)
+x^{-1}p^2 (e_2 \otimes e_2)$
with $x \not= 0;$   or,

  $r= y (e_2 \otimes e_2)$

\end {Lemma}

{ \bf Proof. }
(I)  Let $\{  e_1, e_2, \cdots, e_n \}$  and  $\{  e_1', e_2', \cdots, e_n' \}$
are two basis of $V$  and
    $r = \sum_{i,j= 1 }^{n} k_{ij}(e_{i} \otimes e_{j})$  be strongly symmetric to
the basis   $\{  e_1, e_2, \cdots, e_n \}.$
It is sufficient to show that  $r$ is strongly symmetric to
the basis $\{  e_1', e_2', \cdots, e_n' \}.$
Obviously, there exists $q_{ij} \in k$
such that  $e_{i} = \sum _{s} e_s'q_{si}  $    for $i = 1, 2, \cdots, n.$
By computation,  we have that
$$r = \sum_{ s, t }( \sum _{i,j} k_{ij} q_{si}q_{tj} )(e_{s} \otimes e_{t}).$$
If set   $k_{st}'= \sum _{i,j} k_{ij} q_{si}q_{tj}$,
then
for $l, m, u, v =1, 2, \cdots , n,  $   we have that
 \begin {eqnarray*}
 k_{lm}'k_{uv}'
  &=&  \sum_{i,j,s,t} k_{ij} k_{st} q_{li}q_{mj} q_{us}q_{vt}   \\
  &=&  \sum_{i,j,s,t} k_{is} k_{jt} q_{li}q_{mj} q_{us}q_{vt} { \ \ }
 (  \hbox { \ \  by  }   k_{ij}k_{st} = k_{is}k_{jt}  ) \\
&=&  k_{lu}'k_{mv}'
     \end {eqnarray*}
 and
 \begin {eqnarray*}
 k_{lm}'&=&k_{ml}'
 \end {eqnarray*}
which implies that
 $r$ is strongly symmetric to
the basis $\{  e_1', e_2', \cdots, e_n' \}$

(II)  We only show that  if
 $$ xy =p^2, xz = s^2,  yz=u^2, xu=sp,
p=q, s=t, u=v,$$  then
$$ ys=pu, zp=su.$$
If $p \not= 0,$  then $ysp=yxu=up^2$  and
$usp=uxu = xyz = zp^2,$  which implies
$ys=up$  and $us=zp.$ 
If $p=0,$ then $x=0 $ or $y=0,$ which implies $s=0$ or $u=0.$
                      Consequently,  $ys=up$  and $us=zp.$
$\Box$

\begin {Proposition}\label {1.3}
Let $L$ be a Lie algebra with $dim { \ } L =2$.
Then $r$ is a  solution of CYBE iff
$r$ is strongly symmetric or skew symmetric.
\end {Proposition}

{\bf Proof.}   It is not hard since the tensor $r$ has only 4 coefficients.
 $\Box$

\begin {Proposition}
\label {1.4} Let
$L$ be a Lie algebra with a basis $\{ e_1, e_2, e_3 \}$
such that
$[e_1,e_2]= e_3, [e_2,e_3]=\alpha e_1, [e_3,e_1]=\beta e_2$, where
$\alpha, \beta \in k$. Let $ p, q, s, t, u, v, x, y, z \in k$.

(I)  If $r$ is strongly symmetric or $\alpha, \beta $-skew symmetric
to basis $\{ e_1, e_2, e_3 \},$  then $r$ is a solution of CYBE;

(II)  If  $\alpha \not= 0 , \beta \not=0$,
then $r$ is a   solution of CYBE in $L$
iff   $r$ is strongly symmetric or $\alpha, \beta $-skew symmetric;

(III) If  $\alpha = \beta =0,$    
then $r$ is a solution of CYBE in $L$ iff

 $r= p (e_1 \otimes e_2) +p (e_2 \otimes e_1)+
s (e_1 \otimes e_3)+ t (e_3 \otimes e_1)
+ u (e_2 \otimes e_3)+ v (e_3 \otimes e_2) +
x (e_1 \otimes e_1) +        y (e_2 \otimes e_2)
 +z (e_3 \otimes e_3)$

with  $p \not=0, p^2 = xy, xu = sp, xv =tp, tu =vs,$ or

 $r= s (e_1 \otimes e_3)+ t (e_3 \otimes e_1) + u (e_2 \otimes e_3)+
 v (e_3 \otimes e_2) +
x (e_1 \otimes e_1) + y (e_2 \otimes e_2)
 +z (e_3 \otimes e_3)$

with $xy =xu=xv=ys=yt=0$ and $tu=vs.$

\end {Proposition}

{ \bf Proof .}
Let   $r = \sum_{i,j= 1 }^{3} k_{ij}(e_{i} \otimes e_{j}) \in L \otimes L$
and  $k_{ij} \in k,$ with  $i, j =1, 2, 3.$
It is clear that
\begin {eqnarray*}
\left[ r^{12}, r^{13} \right] &=& \sum _{i,j=1}^3 \sum _{s,t =1}^3 k_{ij}k_{st}
 [e_i,e_s] \otimes e_j \otimes e_t  \\
\left[r^{12}, r^{23} \right] &=& \sum _{i,j=1}^3 \sum _{s,t =1}^3 k_{ij}k_{st}
e_i \otimes[ e_j, e_s ]\otimes e_t \\
\left[r^{13}, r^{23} \right] &=& \sum _{i,j=1}^3 \sum _{s,t =1}^3 k_{ij}k_{st}
e_i \otimes e_s \otimes[e_j, e_t]
\end {eqnarray*}
 By computation, for all $i, j, n =1, 2, 3,$   we have that the coefficient of
  $e_j \otimes e_i \otimes e_i$ in
$[r^{12}, r^{13}] $
is zero and
 $e_i \otimes e_i \otimes e_j$
in  $[r^{13}, r^{23}]$ is zero.

We now see the coefficient
 of  $e_i \otimes e_j \otimes e_n$

   in  $[r^{12}, r^{13}] + [r^{12}, r^{23}] +[r^{13}, r^{23}]$.

(1) $e_1 \otimes e_1 \otimes e_1 (
\alpha k_{12}  k_{31}-
\alpha k_{13} k_{21})$

(2) $e_2 \otimes e_2 \otimes e_2 (
 - \beta k_{21}  k_{32}+\beta k_{23} k_{12} )$

(3) $e_3 \otimes e_3 \otimes e_3 (
k_{31}  k_{23}- k_{32} k_{13})$

(4) $e_1 \otimes e_2 \otimes e_3 (
\alpha k_{22}k_{33} - \alpha k_{32} k_{23} -  \beta k_{11}  k_{33}
+\beta k_{13} k_{13}
+ k_{11}  k_{22} -k_{12} k_{21} )$

(5) $e_2 \otimes e_3 \otimes e_1 (
\beta k_{33}k_{11} - \beta k_{13} k_{31} - k_{22}  k_{11}+ k_{21} k_{21}
+ \alpha k_{22}  k_{33} -\alpha k_{23} k_{32} )$

(6) $e_3 \otimes e_1 \otimes e_2 (
k_{11}k_{22} - k_{21} k_{12} - \alpha k_{33}  k_{22}+ \alpha k_{32} k_{32}
+ \beta k_{33}  k_{11} -\beta k_{31} k_{13} )$

(7) $e_1 \otimes e_3 \otimes e_2 (
-\alpha k_{33}k_{22} + \alpha k_{23} k_{32} + k_{11}  k_{22}- k_{12} k_{12}
-\beta k_{11}  k_{33} + \beta k_{13} k_{31} )$

(8) $e_3 \otimes e_2 \otimes e_1 (
-k_{22}k_{11} + k_{12} k_{21} + \beta k_{33}  k_{11}- \beta k_{31} k_{31}
- \alpha k_{33}  k_{22} + \alpha k_{32} k_{23} )$

(9) $e_2 \otimes e_1 \otimes e_3 (
-\beta k_{11}k_{33} + \beta k_{31} k_{13}  + \alpha k_{22}  k_{33}
- \alpha k_{23} k_{23}
- k_{22}  k_{11} +k_{21} k_{12} )$

(10) $e_1 \otimes e_1 \otimes e_2 (
- \alpha k_{31}k_{22} + \alpha k_{21} k_{32} + \alpha k_{12}  k_{32}
- \alpha k_{13} k_{22})$

(11) $e_2 \otimes e_1 \otimes e_1 (
 -\alpha k_{23}  k_{21}+
\alpha k_{22} k_{31} +\alpha k_{22}  k_{13} - \alpha k_{23} k_{12} )$

(12) $e_1 \otimes e_1 \otimes e_3 (
\alpha k_{21}k_{33} - \alpha k_{31} k_{23} + \alpha k_{12}  k_{33}
 - \alpha k_{13} k_{23} )$

(13) $e_3 \otimes e_1 \otimes e_1 (
 \alpha k_{32}  k_{31}- \alpha k_{33} k_{21}
+ \alpha k_{32}  k_{13} -\alpha k_{33} k_{12} )$

(14) $e_2 \otimes e_2 \otimes e_1 (
- \beta k_{12}k_{31} + \beta k_{32} k_{11} + \beta k_{23}  k_{11}
-\beta  k_{21} k_{31})$

(15) $e_1 \otimes e_2 \otimes e_2 (
 -\beta  k_{11}  k_{32}+
\beta k_{13} k_{12} -\beta k_{11}  k_{23} + \beta k_{13} k_{21} )$

(16) $e_2 \otimes e_2 \otimes e_3 (
\beta k_{32}k_{13} -\beta k_{12} k_{33} -\beta k_{21}  k_{33}
+\beta k_{23} k_{13})$

(17) $e_3 \otimes e_2 \otimes e_2 (
-\beta k_{31}  k_{32}+ \beta k_{33} k_{12}
+ \beta k_{33}  k_{21} - \beta k_{31} k_{23} )$

(18) $e_3 \otimes e_3 \otimes e_1 (
k_{13}k_{21} - k_{23} k_{11} + k_{31}  k_{21}- k_{32} k_{11} )$

(19) $e_3 \otimes e_3 \otimes e_2 (
-k_{23}k_{12} + k_{13} k_{22} - k_{32}  k_{12}+ k_{31} k_{22})$

(20) $e_2 \otimes e_3 \otimes e_3 (
 k_{21}  k_{23} - k_{22} k_{13}
+ k_{21}  k_{32} -k_{22} k_{31} )$

(21) $e_1 \otimes e_3 \otimes e_3 (
- k_{12}  k_{13}+ k_{11} k_{23}
- k_{12}  k_{31} +k_{11} k_{32} )$

(22) $e_1 \otimes e_3 \otimes e_1 (
\alpha k_{23}k_{31} - \alpha k_{33} k_{21} - k_{12}  k_{11}+ k_{11} k_{21}
+\alpha k_{12}  k_{33} -\alpha k_{13} k_{32} )$

(23) $e_1 \otimes e_2 \otimes e_1 (
-\alpha k_{32}k_{21} +\alpha  k_{22} k_{31} - \beta k_{11}  k_{31}
+ \beta k_{13} k_{11} -\alpha k_{13}  k_{22} +\alpha k_{12} k_{23} )$

(24) $e_2 \otimes e_1 \otimes e_2 (
\beta k_{31}k_{12} -\beta k_{11} k_{32} + \alpha k_{22}  k_{32}
-\alpha  k_{23} k_{22} - \beta k_{21}  k_{13} +\beta k_{23} k_{11} )$

(25) $e_2 \otimes e_3 \otimes e_2 (
-\beta k_{13}k_{32} + \beta k_{33} k_{12} + k_{21}  k_{22}- k_{22} k_{12}
-\beta  k_{21}  k_{33} +\beta k_{23} k_{31} )$

(26) $e_3 \otimes e_2 \otimes e_3 (
k_{12}k_{23} - k_{22} k_{13} -\beta k_{31}  k_{33}+ \beta k_{33} k_{13}
+ k_{31}  k_{22} -k_{32} k_{21} )$

(27) $e_3 \otimes e_1 \otimes e_3 (
-k_{21}k_{13} + k_{11} k_{23} + \alpha k_{32}  k_{33}- \alpha k_{33} k_{23}
+ k_{31}  k_{12} -k_{32} k_{11} ).$
    
Let $k_{11}=x, k_{22}=y,  k_{33}=z, k_{12}=p,  k_{21} =q, k_{13}=s,
 k_{31} =t, k_{23}=u,  k_{32}=v.$

It follows from (1)-(27)  that

(28) $\alpha pt= \alpha qs $

(29) $\beta qv= \beta pu$

(30) $tu=vs$

(31) $  \alpha yz - \beta xz + xy - \alpha uv + \beta s^2 - pq =0.$

(32) $ \beta zx -yx+ \alpha yz - \beta st + q^2-\alpha uv =0$

(33) $ xy  -\alpha zy + \beta zx -  pq+ \alpha v^2 - \beta st =0$

(34) $- \alpha zy+ xy- \beta xz+  \alpha uv - p^2+ \beta st =0$

(35) $-xy +\beta xz- \alpha yz + pq - \beta t^2+ \alpha uv =0$

(36) $ -\beta xz+ \alpha yz - yx+  \beta st -  \alpha u^2 + pq =0$

(37)  $\alpha (-ty + qv +pv -sy) = 0$

(38)  $\alpha (-uq + yt +ys -up) = 0$

(39)  $\alpha (qz- tu +pz -su) = 0$

(40)  $\alpha (vt- zq +vs -zp) = 0$

(41)  $\beta (-pt + vx +ux -qt) = 0$

(42)  $\beta (-xv + sp - xu +sq) = 0$

(43)  $\beta (vs-pz +us -qz) = 0$

(44)  $\beta (-tv + zp +zq -tu) = 0$

(45)  $sq -ux + tq- vx = 0$

(46)  $-up + sy -vp +ty = 0$

(47)  $qu-ys + qv  -yt = 0$

(48)  $-ps + xu-pt  +xv = 0.$

(49)  $ \alpha ut - \alpha zq -px +xq +\alpha pz- \alpha sv = 0$

(50)  $- \alpha vq + \alpha yt  -\beta xt + \beta sx - \alpha sy
+\alpha pu = 0$

(51)  $\beta tp - \beta xv + \alpha yv - \alpha uy - \beta qs+ \beta ux = 0$

(52)  $ - \beta sv + \beta zp+ qy - yp - \beta qz+ \beta ut = 0$

(53)  $pu - ys -\beta tz + \beta zs + ty - vq = 0$

(54)  $-qs + xu + \alpha vz -\alpha zu + tp - vx  = 0$
                   
It is clear that $r$ is the solution of CYBE iff relations (28)-(54) hold.

(I)  By computation, we have that if
              $r$ is strongly symmetric or $\alpha, \beta $-skew symmetric
 then  relations (28)-(54) and so  $r$ is a solution of CYBE;

(II)  Let $\alpha \not=0$  and $\beta \not=0.$
 By part (I), we only need show that if $r$ is a solution of  CYBE, then
 $r$ is strongly symmetric or $\alpha, \beta $-skew symmetric.

By computation, we have

(55)  $q^2 =  p^2$   (by  (34)$+$ (32));

(56)  $u^2 =  v^2$  ( by  (33)$+$ (36));

(57)  $t^2 =  s^2$  ( by  (31)$+$ (35) );

(a).If  $s=t \not=0,$ then we have that
 \begin {eqnarray*}
 p&=&q  \hbox { \ \  by (28)} \\
 u&=&v \hbox { \ \ by (30) } \\
 (58) { \ \ \ \ \ \ } yz &=& u^2  \hbox { \ by } (35)-(32) \\
 (59) { \ \ \ \ \ \ } xz &=& s^2  \hbox { \ by } (32)+(33) \\
 (60) { \ \ \ \ \ \ } xy &=& p^2  \hbox { \ by } (36)-(31) \\
(61)  { \ \ \ \ \ \ } xu&=&sp  \hbox { \ by } (41)
 \end {eqnarray*}
then $r$ is strongly symmetric.

 Similarly, we can show that if $p=q\not=0,$ or  $u=v \not=0,$
then $r$  is strongly symmetric.

(b). If  $s = -t \not=0,$
then we have that
   \begin {eqnarray*}
 p&=&-q  \hbox { \ \  by (28)} \\
 u&=&-v \hbox { \ \ by (30) } \\
  y &=& \beta z  \hbox { \ by } (53) \\
 x &=& \alpha z  \hbox { \ by } (50)    \\
 \alpha \beta z^2 + \alpha u^2 + \beta s^2 + p^2 &=& 0   \hbox { \ by } (31)
 \end {eqnarray*}
Then  $r$  is  $\alpha , \beta $-skew symmetric.

Similarly, we have that
if $p=-q \not=0$ or  $u=-v \not=0,$  then
$r$ is $\alpha, \beta $-skew symmetric.

(c). If $s=t=u=v=p=q=0$,
  then we have that
   \begin {eqnarray*}
xz&=&0 \hbox { \ \ by \ }(32)+(33) \\
 yz&=&0 \hbox { \ \ by } (32)+(31)   \\
xy&=& 0 \hbox  { \ \ by } (31)
\end {eqnarray*}
Thus  $r$  is strongly symmetric.

(III) Let $\alpha = \beta =0.$
It is clear that
system of equations  (28)-(54) is equivalent to the below
$$  \left    \{  \begin{array} {l}
 tu=vs, q=p, p^2=xy \\
(s+t)p -(u+v )x =0  \\
(s+t)y-(u+v )p =0  \\
(t-s)y+ (u-v)p  =0  \\
(t-s)p +(u-v )x =0  \\
\end{array} \right.$$

It is easy to check that  if $r$ is  one  case in  part (III)
then the system of equations  hold.
Conversely, if $r$ is a solution of CYBE,
we shall show that $r$ is one of two cases in  part  (III).   
If $p \not=0$ then  $r$  is the first case in part (III).
If $p=0$  then $r$ is the second case in part (III).
        $\Box$

In particular, Proposition \ref {1.4} implies:

\begin {Example} \label {1.5}
 Let $$  sl(2) := \{ x \mid  x  \hbox{ \ is a }  2 \times 2
\hbox { \ matrix
 with trace zero over \  } k   \}$$
 and
 $$e_1 = \left ( \begin {array} {cc}
 0 & 1\\
 1&0
 \end {array}
 \right ),
 e_2 = \left ( \begin {array} {cc}
 0 & -1\\
 1&0
 \end {array}
 \right ),
 e_3 = \left ( \begin {array} {cc}
 2 & 0\\
 0&-2
 \end {array}
 \right ).$$
 Thus  $L$ is a Lie algebra (defined by $[x,y] = xy -yx$)
 with a basis $\{ e_1, e_2, e_3 \}.$ It is clear that
 $$[e_1, e_2] = e_3, [e_2, e_3] = 4e_1, [e_3, e_1] =-4 e_2.$$
Consequently,
 $r$ is a solution of CYBE iff $r$  is strongly symmetric or
4, $-4$-skew symmetric to the basis $\{ e_1, e_2, e_3 \}.$
\end {Example}

\begin {Proposition}
\label {1.6}
$L$ be a Lie algebra with a basis $\{ e_1, e_2, e_3 \}$
such that
$[e_1,e_2]=0,  [e_1,e_3]= e_1+ \beta e_2, [e_2,e_3]=\delta e_2$, where
$\beta, \delta \in k$. Let $ p, q, s, t, u, v, x, y, z \in k$.

(I)  If  $r$ is strongly symmetric,  then
$r$ is a solution of CYBE.

(II)  If  $\beta = 0 , \delta \not=0$,
then $r$   is a solution of CYBE in $L$ iff $r$ is  strongly symmetric, or

  $r=  p(e_1 \otimes e_2)+ q (e_2 \otimes e_1)
+ s (e_1 \otimes e_3)- s (e_3 \otimes e_1)
+u (e_2 \otimes e_3) -u (e_3 \otimes e_2)
+x (e_1 \otimes e_1) +y (e_2 \otimes e_2) $
                                            
with $xu= xs=ys=yu= (1- \delta)us=(1+\delta)s (q +p)
= (1+\delta) u (q +p)= 0$

(III)  If   $\beta \not= 0$ and  $\delta =1,$   then
$r$ is a  solution of CYBE in $L$  iff
 $r$ is strongly symmetric, or

  $r=  p(e_1 \otimes e_2)+ q (e_2 \otimes e_1)
+u (e_2 \otimes e_3) -u (e_3 \otimes e_2)
+x (e_1 \otimes e_1)
+y (e_2 \otimes e_2)$

with $xu= yu=  u (q +p)= 0$

(IV)   If  $ \beta = \delta = 0,$ then
$r$ is a  solution  of CYBE  in $L$ iff
                    
  $r=  p(e_1 \otimes e_2)+ q (e_2 \otimes e_1)
+ s (e_1 \otimes e_3)+s (e_3 \otimes e_1)
+u (e_2 \otimes e_3) +v (e_3 \otimes e_2)
+x (e_1 \otimes e_1)+
y (e_2 \otimes e_2) + z (e_3 \otimes e_3)$
                                            
with $ z \not=0, zp=vs, zq =us, zx =s^2,$  or

  $r=  p(e_1 \otimes e_2)+ q (e_2 \otimes e_1)
+ s (e_1 \otimes e_3)-s (e_3 \otimes e_1)
+u (e_2 \otimes e_3) +v (e_3 \otimes e_2)
+x (e_1 \otimes e_1)+
y (e_2 \otimes e_2)$
                                            
with $ us=vs=xs= xu=xv=0, up=qv,  s(p+q)=0$

\end {Proposition}

{ \bf Proof.}
Let   $r = \sum_{i,j= 1 }^{3} k_{ij}(e_{i} \otimes e_{j}) \in L \otimes L$
and  $k_{ij} \in k,$ with  $i, j =1, 2, 3.$
 By computation, for all $i, j, n = 1, 2 3,$   we have that the coefficient of
  $e_j \otimes e_i \otimes e_i$ in
$[r^{12}, r^{13}] $
is zero and
 $e_i \otimes e_i \otimes e_j$
in  $[r^{13}, r^{23}]$ is zero.

We now see the coefficient
 of  $e_i \otimes e_j \otimes e_n$

   in  $[r^{12}, r^{13}] + [r^{12}, r^{23}] +[r^{13}, r^{23}]$.

(1) $e_1 \otimes e_1 \otimes e_1 (
 -k_{13}k_{11} +  k_{11} k_{31})$;
                                  
(2) $e_2 \otimes e_2 \otimes e_2 (
 -\beta k_{23}k_{12} + \beta k_{21} k_{32}
  - \delta k_{23}k_{22} + \delta k_{22} k_{32})$;

(3) $e_3 \otimes e_3 \otimes e_3 ( 0);$

(4) $e_1 \otimes e_2 \otimes e_3 (
 -k_{32}k_{13} +  k_{12} k_{33} -  \beta k_{13}  k_{13}+  \beta k_{11} k_{33}
- \delta k_{13}  k_{23} + \delta k_{12} k_{33} );$

(5) $e_2 \otimes e_3 \otimes e_1 (
 -\beta k_{33}k_{11} + \beta k_{13} k_{31}
 - \delta k_{33}  k_{21}+ \delta k_{23} k_{31}
 -  k_{23}  k_{31} + k_{21} k_{33} );$

(6) $e_3 \otimes e_1 \otimes e_2 (
 -k_{33}k_{12} +  k_{31} k_{32} -  \beta k_{33}  k_{11}+  \beta k_{31} k_{13}
- \delta k_{33}  k_{12} + \delta k_{32} k_{13} );$

(7) $e_1 \otimes e_3 \otimes e_2 (
 -k_{33}k_{12} +  k_{13} k_{32} -  \beta k_{13}  k_{32}+  \beta k_{11} k_{33}
- \delta k_{13}  k_{32} + \delta k_{12} k_{33} );$

(8) $e_3 \otimes e_2 \otimes e_1 (
 -\beta k_{33}k_{11} + \beta k_{31} k_{31} -  \delta k_{33}  k_{21}
 +  \delta k_{32} k_{31}
- k_{33}  k_{21} +  k_{31} k_{23} );$

(9) $e_2 \otimes e_1 \otimes e_3 (
 -\beta k_{31}k_{13} + \beta  k_{11} k_{33}
 - \delta k_{31}  k_{23}+  \delta k_{21} k_{33}
-  k_{23}  k_{13} +  k_{21} k_{33} );$

(10) $e_1 \otimes e_1 \otimes e_2 (
 -k_{31}k_{12} +  k_{11} k_{32} -   k_{13}  k_{12}+   k_{11} k_{32});$

(11) $e_2 \otimes e_1 \otimes e_1 (
 -k_{23}k_{11} +  k_{21} k_{31} -   k_{23}  k_{11}+  k_{21} k_{13});$

(12) $e_1 \otimes e_1 \otimes e_3 (
 -k_{31}k_{13} +  k_{11} k_{33} -   k_{13}  k_{13}+   k_{11} k_{33}
 );$

(13) $e_3 \otimes e_1 \otimes e_1 (
 -k_{33}k_{11} +  k_{31} k_{31} -   k_{33}  k_{11}+  k_{31} k_{13}
);$

(14) $e_2 \otimes e_2 \otimes e_1 (
 -\beta k_{32}k_{11} +  \beta k_{12} k_{31} -  \delta k_{32}  k_{21}
 +  \delta k_{22} k_{31}
- \beta k_{23}  k_{11} + \beta k_{21} k_{31}
- \delta k_{23}  k_{21} + \delta k_{22} k_{31}
 );$

(15) $e_1 \otimes e_2 \otimes e_2 (
 -\beta k_{13}k_{12} + \beta  k_{11} k_{32} -  \delta k_{13}  k_{22}
 +  \delta k_{12} k_{32}
- \beta k_{13}  k_{21} + \beta k_{11} k_{23}
- \delta k_{13}  k_{22} + \delta k_{12} k_{23}
 );$

(16) $e_2 \otimes e_2 \otimes e_3 (
 -\beta k_{32}k_{13} +  \beta k_{12} k_{33} -  \delta k_{32}  k_{23}
 + \delta k_{22} k_{33}
- \beta k_{23}  k_{13} + \beta k_{21} k_{33}
- \delta k_{23}  k_{23} + \delta k_{22} k_{33}
);$

(17) $e_3 \otimes e_2 \otimes e_2 (
 -\beta k_{33}k_{12} + \beta  k_{31} k_{32}- \delta k_{33}k_{22}
 +\delta k_{32}k_{32}
  -  \beta k_{33} k_{21} + \beta k_{31}  k_{23}
  - \delta k_{33} k_{22}
  + \delta k_{32} k_{23}
  );$

(18) $e_3 \otimes e_3 \otimes e_1 (0);$

(19) $e_1 \otimes e_3 \otimes e_3 (0);$

(20) $e_3 \otimes e_3 \otimes e_2 (0);$

(21) $e_2 \otimes e_3 \otimes e_3 (0);$

(22) $e_1 \otimes e_3 \otimes e_1 (
 -k_{33}k_{11} +  k_{13} k_{31} -   k_{13}  k_{31}+   k_{11} k_{33})
 = e_1 \otimes e_3 \otimes e_1 (0)$;

(23) $e_1 \otimes e_2 \otimes e_1 (
 -k_{32}k_{11} +  k_{12} k_{31} -  \beta k_{13}  k_{11}+  \beta k_{11} k_{31}
- \delta k_{13}  k_{21} + \delta k_{12} k_{31}
-  k_{13}  k_{21} +  k_{11} k_{23});$

(24) $e_2 \otimes e_1 \otimes e_2 (
 -\beta k_{31}k_{12} + \beta k_{11} k_{32}
 - \delta k_{31}  k_{22}+  \delta k_{21} k_{32}
- k_{23}  k_{12} +  k_{21} k_{32}
 -\beta k_{23}k_{11} + \beta k_{21} k_{13}
  - \delta k_{23}  k_{12}+  \delta k_{22} k_{13}  );$

(25) $e_2 \otimes e_3 \otimes e_2 (
 -\beta k_{33}k_{12} + \beta k_{13} k_{32}
  -  \delta k_{33}  k_{22}+  \delta k_{23} k_{32}
- \beta k_{23}  k_{31} + \beta k_{21} k_{33}
 + \delta k_{33}  k_{22}-  \delta k_{22} k_{33}
              );$

(26) $e_3 \otimes e_2 \otimes e_3 (
 -\beta k_{33}k_{13} + \beta k_{31} k_{33}
  - \delta k_{33}  k_{23}+  \delta k_{32} k_{33}
);$

(27) $e_3 \otimes e_1 \otimes e_3 (
 - k_{33}k_{13} +   k_{31} k_{33}
 );$

Let $k_{11}=x, k_{22}=y,  k_{33}=z, k_{12}=p,  k_{21} =q, k_{13}=s,
 k_{31} =t, k_{23}=u,  k_{32}=v.$

   It follows from (1)-(27) that

(28) $-sx + xt=0$;

(29) $-\beta up +\beta qv -\delta uy + \delta yv=0$;

(30) $- vs + pz -\beta s^2 + \beta xz -\delta su + \delta zp =0$;

(31) $-\beta xz + \beta st - \delta zq + \delta ut - ut   + qz =0$;

(32) $- zp + tv -\beta zx + \beta ts -\delta zp + \delta vs =0$;

(33) $- zp + sv -\beta st + \beta xz -\delta sv + \delta zp =0$;

(34) $- \beta zx + \beta tt - \delta zq + \delta vt - zq + tu =0$;

(35) $- \beta st + \beta xz  -\delta tu + \delta qz -us + qz=0 $;

(36) $- tp + xv - sp + vx =0$;

(37) $- ux + qt-ux + qs =0$;

(38) $- st + xz - s^2 +  xz =0$;

(39) $-zx +tt -zx +st =0$;

(40) $-\beta vx +\beta pt -\delta vq + \delta yt -\beta ux + \beta qt
-\delta uq + \delta yt =0$;

(41) $-\beta sp + \beta xv  -\delta  sy + \delta pv -\beta sq + \beta xu
- \delta sy + \delta pu =0$;

(42) $-\beta vs + \beta pz -\delta vu + \delta yz -\beta su + \beta zq
-\delta u^2 + \delta yz  =0$;

(43) $-\beta  zp + \beta tv - \delta zy + \delta vv  -\beta zq + \beta tu
-\delta zy +\delta vu =0$;

(44) $- vx + pt -\beta sx + \beta xt -\delta sq + \delta pt
-sq + xu =0$;

(45) $- \beta pt + \beta vx -\delta ty + \delta qv  -up + qv
-\beta ux + \beta qs - \delta up + \delta ys =0$;

(46) $-\beta zp + \beta sv -\beta ut + \beta qz =0$;

(47) $- \beta zs + \beta tz -\delta zu + \delta vz =0$;

(48) $- zs + tz =0$;

It is clear that   $r$   is a solution of CYBE iff
 (28)-(48) hold.

(I)  It is trivial.

(II)  Let  $ \beta =0$  and  $\delta \not=0.$
We only show that if  $r$ is a solution of CYBE, then $r$
is one of two cases in part (II).
If $z\not=0$, we have that
$u = v, s=t, yz = u^2, xz = s^2$ and $ux = sp $   by (47), (48), (42),
(38) and (36),
respectively. It follows from (30) and (33) that $pz = us$ and
from (31) and (34) that $qz = us.$  Thus $q=p$ and $xy =p^2$
since $xyz =xu^2=spu =zp^2. $  By Lemma \ref {1.2}, $r$ is strongly
symmetric.

If $z=0$, we have that $u=-v$ and $s=-t$ by (43), (42), (38) and (39). It
follows from (28), (29), (36), (40), (30), (44) and (45) that
$sx =0, yu=0, xv=0, ys=0, (1-\delta)us=0,
(1+\delta)s(p+q)=0$ and $(1+\delta)u(p+q)=0,$  respectively. Consequently,
$r$ is the second case.

       (III) Let  $\beta \not=0,$ $\delta =1.$  
We only show that if $r$  is a solution of CYBE, then  $r$ is one case in part
(III).
If $z\not=0$, we have that
$t=s , u=v, q=p, xz =s^2, zp = su, sp =xu$ and $u ^2 = yz$ by
(48), (47), (46), (38), (32), (36) and (42), respectively. Since
$xyz =xu^2 = spu = zp^2$, $xy=p^2.$   Thus
 $r$ is strongly symmetric by Lemma \ref {1.2}.

If $z=0$, then $st =0, t=0$ and $s=0$ by (31), (39) and (38), respectively.
It follows from (42) and (43) that $u=-v$. By (36), (45) and (29), $xu=0,
u(p+q)=0$ and $uy=0$. Consequently $r$  is the second case.

(IV)        Let $\beta = \delta =0.$
We only show that if  $r$ is a solution of CYBE, then
$r$ is one case in part  (IV).
If $z\not=0$, we have that $s=t, s^2=xz, zp=vs $ and $zq=us$
by (48), (38), (32) and (34), respectively. Consequently
 $r$ is the first case in part (IV). If $z=0$,
 then $s= -t$ by (38) and (39). It follows from (30), (31), (28), (37), (36),
 (45) and (44) that $vs=0, us=0, xs=0, xu=0, xv=0 , up = qv$ and $s(p+q)=0.$
Consequently,  $r$
is the second case in part (IV).
$\Box$

\begin {Corollary}\label {1.7}
Let $L$ is a Lie algebra with $dim { \ } L \le 3$  and $r \in L \otimes L.$
If $r$ is strongly symmetric, then $r$ is a solution of CYBE in $L$.
                                     
\end{Corollary}

{\bf Proof .}   If $k$ is algebraically closed, then  $r$ is a solution of
 CYBE by Proposition  \ref {1.4}, \ref {1.6} and
 \ref {1.3}.
 If $k$ is not algebraically closed, let
  $P$ be algebraically closure of $k$.
  We can construct a Lie algebra   $ L_P = P \otimes L$ over $P$,
  as in \cite [section 8] {J}.
Set $ \Psi : L \longrightarrow  L_P $  by sending $x $  to
  $1 \otimes x$. It is clear that $L_P$ is a Lie algebra over
  $P$ and  $\Psi$  is  homomorphic with $ker \Psi =0$ over $k.$
  Let  $$\bar r = (\Psi \otimes \Psi )(r).$$
  Obviously, $\bar r$  is  strongly symmetric. Therefore  $\bar r$
 is a solution of CYBE in $L_P$ and so is $r$  in $L$.
$\Box$

\section {Coboundary Lie bialgebras}

In this section, using the general solution, which are obtaied in the section
above,  of CYBE in Lie algebra $L$ with $dim \ \ L \le 3, $ we give the
       the sufficient and  necessary conditions which  $(L, \hbox {[ \ ]},
 \Delta _r, r)$ is a coboundary
 (or triangular ) Lie bialgebra.

We now observe the  connection between solutions
of CYBE and triangular Lie bialgebra structures. It is clear that if
 $(L, [ \hbox { \ } ], \Delta _r, r)$   is a triangular
Lie bialgebra then  $r$ is a solution of CYBE.
Conversely, if $r$ is a solution of CYBE and $r$ is skew symmetric with
$r \in L \otimes L$ , then
               $(L, [ \hbox { \ } ], \Delta _r, r)$   is a triangular
Lie bialgebra  by \cite [Proposition 1] {T} .
\begin {Theorem}
\label {2.1}
Let $L$ be a  Lie algebra with
       a basis $\{ e_1, e_2, e_3 \}$
such that
$[e_1,e_2]= e_3, [e_2,e_3]=\alpha e_1, [e_3,e_1]=\beta e_2$, where
$\alpha, \beta \in k$ and $\alpha \beta \not=0$ or $\alpha =\beta =0$.
Let $p, u, s \in k$ and
      $r=  p(e_1 \otimes e_2)- p (e_2 \otimes e_1)
+ s (e_1 \otimes e_3)- s (e_3 \otimes e_1)
+ u (e_2 \otimes e_3)- u (e_3 \otimes e_2)$.
Then

(I)   $(L, [\hbox { \ }], \Delta _r, r )$ is  a coboundary Lie bialgebra
iff $r$  is skew symmetric;

(II)       $(L, [\hbox { \ }],\Delta _r, r )$ is  a
triangular  Lie bialgebra
iff  $$\beta s^2 + \alpha u^2 + p^2 =0, $$
\end {Theorem}

 {\bf Proof.}
(I)
It is sufficient to show that
 $$ (1+ \xi + \xi ^2) (1 \otimes \Delta ) \Delta (e_i) =0$$
for $i=1, 2, 3.$
First, by computation, we have that
\begin {eqnarray*}
&{ \ }& (1 \otimes \Delta ) \Delta (e_1) \\
&=&
\{ ( e_1 \otimes e_2\otimes e_3 ) (\beta ps -  \beta ps)
 +(e_3 \otimes e_1\otimes e_2)(\beta ps)
 + (e_2 \otimes e_3\otimes e_1) (- \beta sp) \} \\
&+&\{ ( e_1 \otimes e_3\otimes e_2 ) (-\beta ps + \beta ps)
 +(e_3 \otimes e_2\otimes e_1)( - \beta ps)
 + (e_2 \otimes e_1\otimes e_3) (\beta sp) \} \\
&+&\{ ( e_1 \otimes e_1\otimes e_2 ) (\alpha \beta su )
 +(e_1 \otimes e_2\otimes e_1)( - \beta \alpha su)
 + (e_2 \otimes e_1\otimes e_1) (0) \} \\
&+& \{ ( e_1 \otimes e_1\otimes e_3 ) ( - \alpha pu)
 +(e_1 \otimes e_3\otimes e_1)(\alpha  pu)
 + (e_3 \otimes e_1\otimes e_1) (0) \} \\
&+&\{   ( e_2 \otimes e_1\otimes e_2 ) (- \beta ^2 s^2)
 +(e_1 \otimes e_2\otimes e_2)(0)
 + (e_2 \otimes e_2\otimes e_1) ( \beta ^2 s^2) \} \\
&+&\{  ( e_3 \otimes e_3\otimes e_1 ) ( p^2 )
 +(e_3 \otimes e_1\otimes e_3)(- p^2)
 + (e_1 \otimes e_3\otimes e_3) (0) \}
\end {eqnarray*}
Thus
 $$ (1+ \xi + \xi ^2) (1 \otimes \Delta ) \Delta (e_1) =0$$
Similarly, we have that
 $$ (1+ \xi + \xi ^2) (1 \otimes \Delta ) \Delta (e_2) =0$$
           $$ (1+ \xi + \xi ^2) (1 \otimes \Delta ) \Delta (e_3) =0$$
Thus  $(L, [\hbox { \ }],\Delta _r,r)$ is a coboundary Lie bialgebra.

(II) By part (I), it is sufficient to show that $r$ is a solution of CYBE
iff
  $$\beta s^2 + \alpha u^2 + p^2 =0. $$
(a). Let $\alpha \not=0 $ and $\beta \not=0.$ Considering $r$ is skew
symmetric, by Proposition \ref {1.4} (II), we have
                             $r$ is a solution of CYBE
iff
  $$\beta s^2 + \alpha u^2 + p^2 =0. $$
(b). Let $\alpha = \beta=0.$                  Considering $r$ is skew
symmetric, by Proposition \ref {1.4} (III), we have
                              $r$ is a solution of CYBE
iff
  $$\beta s^2 + \alpha u^2 + p^2 =0. $$
$\Box$

Let $(L,[ \hbox{ \ }])$  be a familiar Lie algebra, namely,
 Euclidean 3-space
under vector cross product.
By Theorem \ref {2.1},  $(L, [ \hbox { \ }], \Delta _r, r) $
 is not  a triangular Lie bialgebra for any   $0 \not= r \in L \otimes L.$
In fact, \cite [Example 2.14] {M} already  contains the essence of this
observation.

\begin {Example} \label {2.2} Under Example \ref {1.5}, we have the
following:
                                                                         
(i)   $(sl(2), [\hbox { \ }],\Delta _r, r )$ is  a coboundary Lie bialgebra
iff $r$ is skew symmetric;

(ii)       $(sl(2), [\hbox { \ }],\Delta _r, r )$ is  a
triangular  Lie bialgebra iff
$$-4 s^2 +4 u^2 + p^2 =0, $$         
where    $r=  p(e_1 \otimes e_2)- p (e_2 \otimes e_1)
+ s (e_1 \otimes e_3)- s (e_3 \otimes e_1)
+u (e_2 \otimes e_3) -u (e_3 \otimes e_2)$  and $p, s, u \in k,$
\end {Example}

 \begin {Theorem} \label  {2.3}
Let  $L$ be a Lie algebra
 with a basis   $\{ e_1, e_2, e_3 \}$  such that
$$[e_1e_2]= 0, [e_1, e_3]= e_1 + \beta e_2, [e_2, e_3]=\delta e_2,$$
where $\delta, \beta \in k.$  Let $p, s, u \in k$ and
 $r=  p(e_1 \otimes e_2)- p (e_2 \otimes e_1)
+ s (e_1 \otimes e_3)- s (e_3 \otimes e_1)
+ u (e_2 \otimes e_3)- u (e_3 \otimes e_2)$

(I)  $(L, [\mbox { \ }],\Delta _r, r )$ is a coboundary Lie bialgebra
iff $$ (\delta +1) ((\delta -1)u + \beta s)s =0; $$

(II) If  $\beta =0,$  then
 $(L, [\mbox { \ }],\Delta _r, r )$ is a triangular Lie bialgebra
 iff $$(1- \delta )  us =0.$$

(III) If   $\beta \not=0$ and  $\delta =1,$ then
 $(L, [\mbox { \ }],\Delta _r,r)$ is a triangular Lie bialgebra iff
 $$ s =0.$$
iff   $(L, [\mbox { \ }],\Delta _r, r )$ is a coboundary Lie bialgebra.
\end {Theorem}
{ \bf Proof }.
 (I) We get by computation
 \begin {eqnarray*}
&{ \ }& ( 1 \otimes \Delta ) \Delta (e_1) \\
 &=& \{ ( e_1 \otimes e_1\otimes e_2 ) (\beta \delta ss -  \delta us)
 +(e_1 \otimes e_2\otimes e_1)(- \beta \delta ss + \delta us)
 + (e_2 \otimes e_1\otimes e_1) (0) \} \\
&+&  \{ ( e_2 \otimes e_2\otimes e_1 ) (-\beta us + uu+ \beta ^2 s^2
- \beta su)                                       \\
 &+& (e_2 \otimes e_1\otimes e_2)( \beta us -uu + \beta us - \beta ^2 s^2)
  + (e_1 \otimes e_2\otimes e_2) (0) \} \\
&{ \ }& ( 1 \otimes \Delta ) \Delta (e_2) \\
 &=& \{ ( e_1 \otimes e_1\otimes e_2 ) ( \delta ^2 ss)
 +(e_1 \otimes e_2\otimes e_1)(- \delta ^2 s^2 )
 + (e_2 \otimes e_1\otimes e_1) (0) \} \\
 &+& \{ ( e_2 \otimes e_2\otimes e_1 ) ( \delta \beta ss - \delta su)
  +(e_2 \otimes e_1\otimes e_2)( \delta  su - \delta \beta ss )
 + (e_1 \otimes e_2\otimes e_2) (0) \} \\  
&{ \ }& ( 1 \otimes \Delta ) \Delta (e_3)  \\
 &=& \{ ( e_1 \otimes e_2\otimes e_3 ) (\delta su +  \beta ss)
 +(e_2 \otimes e_3\otimes e_1)(-\delta us - \beta ss) \\
 &+& (e_3 \otimes e_1\otimes e_2) ( \delta \delta us + \beta \delta ss
 + \beta ss -su) \} \\
&+&\{ ( e_1 \otimes e_3\otimes e_2 ) (-\delta su-  \beta ss)
 +(e_3 \otimes e_2\otimes e_1)( - \beta \delta  ss - \delta \delta us
 -\beta ss + su ) \\
 &+&  (e_2 \otimes e_1\otimes e_3) (\delta us + \beta ss) \} \\
&+& \{ ( e_1 \otimes e_1\otimes e_2 ) (- \delta  sp - \delta \delta ps
+ sp + \delta sp ) \\
 &+& (e_1 \otimes e_2\otimes e_1)( \delta ps + \delta \delta ps -\delta sp
 -sp) +
 (e_2 \otimes e_1\otimes e_1) (0) \} \\
&+& \{
 ( e_1 \otimes e_1\otimes e_3 ) ( ss)
 +(e_1 \otimes e_3\otimes e_1)(-ss)
 + (e_3 \otimes e_1\otimes e_1) (0) \} \\
&+& \{   ( e_2 \otimes e_2\otimes e_1 ) (- \beta  ps + pu - \delta \beta sp
-\delta \delta up - \delta \beta ps + \delta pu - \beta sp - \delta up)  \\
&+& (e_2 \otimes e_1\otimes e_2)( \beta  ps - pu + \delta \beta sp
+\delta \delta up + \delta \beta ps - \delta pu + \beta sp + \delta up) ) \\
&-& (e_1 \otimes e_2\otimes e_2) (0) \} \\
&+& \{ ( e_2 \otimes e_2\otimes e_3 ) ( \delta \delta uu + \delta \beta us
+ \beta \delta us + \beta \beta ss )   \\
 &+& (e_2 \otimes e_3\otimes e_2)(   -\delta \delta uu - \delta \beta us
- \beta \delta us - \beta \beta ss )
                    +        (e_3 \otimes e_2\otimes e_2) (0) \}
\end {eqnarray*}
Consequently,   $$ (1+ \xi + \xi ^2) (1 \otimes \Delta ) \Delta (e_1) =0$$
 $$ (1+ \xi + \xi ^2) (1 \otimes \Delta ) \Delta (e_2) =0$$
 and
 $$ (1+ \xi + \xi ^2) (1 \otimes \Delta ) \Delta (e_3) =0$$
iff
  $$ \delta ^2 us + \delta \beta s^2 + \beta s^2 - us =0.$$
This implies that
    $(L, [\mbox { \ }],\Delta _r,r)$ is a coboundary Lie bialgebra
    iff
 $$ \delta ^2 us + \delta \beta s^2 + \beta s^2 - us =0   { \ \ \ \ }$$

 (II) If $(L, [ \hbox { \ } ], \Delta _r, r)$  is a triangular Lie bialgebra,
  then
 $(\delta ^2 -1)us =0$ by part (I). If $(\delta +1) \not=0,$ then
 $( 1 - \delta )us =0.$  If $\delta +1 =0, $  then $(1 - \delta )us =0$
 by Proposition \ref {1.6} (II). Conversely,
 if $(1 - \delta )us =0$, then we have
 that $(L, [ \hbox { \ } ], \Delta _r, r)$ is a coboundary Lie bialgebra
 by part (I).
 Since $r$ is skew symmetric we have that $r$ is a solution of CYBE by
 Proposition   \ref {1.6} (II) (IV). Thus
 $(L, [ \hbox { \ } ], \Delta _r, r)$ is a triangular Lie bialgebra.

(III)  Similarly, we can show that part (III) holds by part (I) and
Proposition \ref {1.6} (III). $\Box$

\begin {Theorem}  \label   {2.4}
If $L$ is a Lie algebra with $dim L =2$ and $r \in L \otimes L$, then
  $(L, [\mbox { \ }],\Delta _r,r)$ is a triangular  Lie bialgebra   iff
  $(L, [\mbox { \ }],\Delta _r,r)$ is a coboundary   Lie bialgebra
  iff   $r$ is skew symmetric.
 \end {Theorem}
 { Proof.}  It is an immediate consequence of the main result of \cite {M}.
  $\Box$

\vskip 2cm

{\bf Acknowledgement }  I would like to show my gratitude to Professor
Yonghua Xu and  referee
  for their guidance and help.

\end {document}